# Encodings of trajectories and invariant measures

G. S. Osipenko

**Abstract.** We consider a discrete dynamical system on a compact manifold $M$ generated by a homeomorphism $f$. Let $C = \{M(i)\}$ be a finite covering of $M$ by closed cells. The symbolic image of a dynamical system is a directed graph $G$ with vertices corresponding to cells in which vertices $i$ and $j$ are joined by an arc $i \to j$ if the image $f(M(i))$ intersects $M(j)$. We show that the set of paths of the symbolic image converges to the set of trajectories of the system in the Tychonoff topology as the diameter of the covering tends to zero. For a cycle on $G$ going through different vertices, a simple flow is by definition a uniform distribution on arcs of this cycle. We show that simple flows converge to ergodic measures in the weak topology as the diameter of the covering tends to zero.

Bibliography: 28 titles.

**Keywords:** pseudotrajectory, recurrent trajectory, chain recurrent set, ergodic measure, symbolic image, flow on a graph.

## § 1. Introduction

Let $f\colon M \to M$ be a homeomorphism of a compact Riemannian manifold $M$, which generates the discrete dynamical system

$$x_{n+1} = f(x_n); \tag{1.1}$$

let $\rho(x,y)$ be the distance on $M$.

The basic concept underlying the analysis which follows is the concept of the symbolic image of a dynamical system (see [1] and [2]). This concept combines symbolic dynamics (see [3]–[8]) and numerical methods (see [9]).

Let $C = \{M(1), \dots, M(n)\}$ be a finite covering of the manifold $M$ by closed sets; the set $M(i)$ will be referred to as the cell with index $i$.

**Definition 1** (see [10]). The *symbolic image of the dynamical system* (1.1) for a covering $C$ is a directed graph $G$ with vertices $\{i\}$ corresponding to the cells $\{M(i)\}$. Two vertices $i$ and $j$ are joined by a directed edge (an arc) $i \to j$ if and only if

$$f(M(i)) \cap M(j) \neq \varnothing.$$

This research was carried out with the financial support of the Russian Foundation for Basic Research (grant no. 19-01-00388-a).

*AMS 2020 Mathematics Subject Classification.* Primary 37C50.

The symbolic image generates the symbolic dynamics, which reflects the dynamics of system (1.1); using it we can obtain useful information about the global structure of the dynamics of the system. The symbolic image depends on the covering $C$, if it is changed the symbolic image changes. The existence of the edge $i \to j$ guarantees that the cell $M(i)$ contains a point $x$ whose image $f(x)$ lies in $M(j)$. In other words, the edge $i \to j$ is the trace of the mapping $x \to f(x)$, where $x \in M(i)$ and $f(x) \in M(j)$. If there is no such edge $i \to j$, then there is no point $x \in M(i)$ whose image $f(x)$ lies in $M(j)$.

We shall consider coverings $C$ with polyhedral cells $M(i)$ that intersect in boundary discs. Such coverings always exist — this follows from the theorem on the triangulation of a compact manifold (see [11]). In numerical calculations (see [1]), $M$ is a compact domain in $R^d$ and cells $M(i)$ are cubes or parallelepipeds. Let $d = \mathrm{diam}(C)$ be the greatest of the diameters of the cells of the covering $C$. It is called the diameter of the covering $C$. In sections pertaining to measure theory, we shall change from the covering $C$ to a partition $C^*$ by attributing each boundary disc to one of the adjoining cells. In this case $C^*$ is a measurable partition of the manifold $M$.

Recall that a doubly infinite sequence of points $T = \{x_n\}$ is a trajectory of system (1.1) if $f(x_n) = x_{n+1}$. A doubly infinite sequence of points $\chi = \{x_n\}$ is called an $\varepsilon$-trajectory or a pseudotrajectory if $\rho(f(x_n), x_{n+1}) < \varepsilon$ for any $n$.

**Definition 2.** A doubly infinite sequence $\omega = \{z_k\}$ of vertices of a graph $G$ is called a *path* (or an *admissible path*) if, for each $k$, the graph $G$ contains the arc $z_k \to z_{k+1}$.

We denote the vertex set of the graph $G$ by $V$. The symbolic image of $G$ can be considered as a set-valued mapping between vertices, $G\colon V \to V$, where the image $G(i)$ is the set of vertices $j$ which are the terminal points of arcs $i \to j$:

$$G(i) = \{j\colon i \to j\}.$$

There exists a natural set-valued mapping $h\colon M \to V$ from $M$ onto the vertex set $V$ of the symbolic image which associates with each point $x$ the family of vertices $i$ such that $x \in M(i)$; that is,

$$h(x) = \{i\colon x \in M(i)\}.$$

By the definition of the symbolic image, the diagram

$$\begin{array}{ccc} M & \xrightarrow{f} & M \\ \downarrow h & & \downarrow h \\ V & \xrightarrow{G} & V \end{array} \qquad (1.2)$$

is commutative in the following sense:

$$h(f(x)) \subset G(h(x)). \qquad (1.3)$$

In fact, let $i \in h(x)$ and $j \in h(f(x))$. Then $M(j) \cap f(M(i)) \neq \varnothing$ and there exists an arc $i \to j$, which means that $j \in G(i)$ or $h(f(x)) \in G(h(x))$. Therefore, $h(f(x)) \subset G(h(x))$. Note that the equality $h(f(x)) = G(h(x))$ cannot be guaranteed. However, the inclusion (1.3) is sufficient for $h$ to transform trajectories of (1.1) into admissible paths of the symbolic image,

$$h(T) = \omega = \{i_n \colon f^n(x) \in M(i_n)\}.$$

In this case we say that the path $\omega$ is the trace of the trajectory $T$ on the symbolic image $G$. The trace $\omega$ can be interpreted as an encoding of the trajectory $T$. Let $P$ be the set of admissible paths and let Cod be the set of encodings of trajectories (note that Cod $\subset P$).

If the symbolic image $G$ contains a path $\omega = \{i_n\}$, then we can define the sequence of points $\chi = \{x_n \colon x_n \in M(i_n)\}$, which is a pseudotrajectory. In this case the pseudotrajectory $\chi$ is said to be the trace of the path $\omega$. It is clear that the trace of a path is not uniquely defined. The first result in our paper (Theorem 2) gives conditions for a trace $\chi$ to converge to some trajectory $T$, provided that the diameter of the covering converges to zero. We will consider the case of convergence to a recurrent trajectory.

Let $\Lambda$ be the set of traces of admissible paths on $G$ and let Tr be the set of trajectories, Tr $\subset \Lambda$. One result in our paper (Theorem 5) asserts that the set of traces of paths of the symbolic image converges to the set of trajectories of the system in the Tychonoff topology as the diameter of the covering tends to zero.

A measure $\mu$ on $M$ is said to be invariant for $f$ (or $f$-invariant) if

$$\mu(f^{-1}(A)) = \mu(A) = \mu(f(A))$$

for any measurable set $A \subset M$. The Krylov-Bogolyubov theorem (see [12] and [13]) guarantees the existence of an invariant measure $\mu$ which is normalized on $M$: $\mu(M) = 1$. We denote the set of all $f$-invariant normalized measures by $\mathscr{M}(f)$. If $\mu_1, \mu_2$ lie in $\mathscr{M}(f)$ and if $\alpha, \beta \geqslant 0$ and $\alpha + \beta = 1$, then $\mu = \alpha\mu_1 + \beta\mu_2$ also lies in $\mathscr{M}(f)$ (see [13]). In this case we say that $\mu$ is a sum of the measures $\mu_1$ and $\mu_2$. The set of all $f$-invariant measures $\mathscr{M}(f)$ is a convex weakly compact set; see [14]. The weak convergence $\mu_n \to \mu$ means that

$$\int_M \varphi \, d\mu_n \to \int_M \varphi \, d\mu$$

for any continuous function $\varphi \colon M \to R$.

The measure $\mu$ is ergodic if either $\mu(A) = 0$ or $\mu(M \setminus A) = 0$ for any invariant measurable set $A$. A point $\mu$ in a convex compact set $\mathscr{M}$ is an extreme point of $\mathscr{M}$ if $\mu$ cannot be represented as $\mu = \alpha\mu_1 + (1-\alpha)\mu_2$, where $\mu_1, \mu_2 \in \mathscr{M}$, $\mu_1 \neq \mu_2$ and $0 < \alpha < 1$. We let ext$(\mathscr{M})$ denote the set of extreme points of a convex compact set $\mathscr{M}$. In [15] it was shown that the ergodic measures are exactly the extreme points of the convex set $\mathscr{M}(f)$. The set of ergodic measures also has other interesting properties in the weak topology; see[16]–[18].

**Definition 3.** Let $G$ be an directed graph. A distribution $m = \{m_{ij}\}$ on the edges $\{i \to j\}$ is called a *flow* if:
  1) $m_{ij} \geqslant 0$;
  2) $\sum m_{ij} = 1$;
  3) $\sum_k m_{ki} = \sum_j m_{ij}$ for each $i \in G$.

The equality in 3) can be interpreted as Kirchhoff's law: the flow coming into the vertex $i$ is equal to the flow going out of $i$. Each invariant measure $\mu$ generates a flow on the symbolic image by the formula

$$m_{ij} = \mu(f(M(i)) \cap M(j)) = \mu(f^{-1}(M(j)) \cap M(i));$$

for details, see [19]. This flow is called the trace of the measure $\mu$ on the symbolic image $G$. If there is a flow $m$ on the symbolic image $G$, then a measure $\mu$ on $M$ can be defined by

$$\mu(A) = \sum_i m_i \frac{v(A \cap M(i))}{v(M(i))}, \qquad (1.4)$$

where $A$ is a measurable set, $m_i = \sum_j m_{ij}$ is the measure of the cell $i$, and $v(\cdot)$ is the Lebesgue measure. Here we assume that $v(M(i)) \neq 0$. The measure $\mu(M(i))$ of the cell $M(i)$ coincides with $m_i$. The measure $\mu$ thus constructed is called the trace of the flow $m$. This measure is not invariant, but according to [19] it approximates an invariant measure as the diameter of the covering converges to zero. Let $\Omega(G)$ be the set of traces of flows on $G$. Theorem 5 in §2 is the analogue of Theorem 4 in [19] for invariant measures: $\Omega(G)$ converges to $\mathscr{M}(f)$ in the weak topology as the diameter of the covering tends to zero.

The set $\mathscr{M}(G)$ of flows on the graph $G$ has the natural structure of a convex set. Consider a cycle $\omega = \{i_1 \to i_2 \to \cdots \to i_p \to i_1\}$ on $G$; we say that $\omega$ is a simple cycle if all the vertices $i_n$, $n = 1, 2, \ldots, p$, are distinct. A simple cycle generates a simple flow $m(\omega)$ localized on $\omega$:

$$m_{ij} = \frac{1}{p} \quad \text{for } (i \to j) \in \omega, \qquad m_{ij} = 0 \quad \text{for } (i \to j) \notin \omega.$$

In [19] it was shown that any flow can be represented as a sum of simple ones; that is, the simple flows form the set of extreme points of $\mathscr{M}(G)$. We shall show (see Theorem 9 below) that the traces of simple flows converge to ergodic measures as the diameter of the covering tends to zero.

## §2. Shadowing of trajectories

Let $C$ be a covering, let $G$ be a symbolic image, let $\operatorname{diam} f(M(i)) = q_i$ be the diameter of the image of a cell, and let $q = \max q_i$ be the greatest of the diameters of the images of cells. If $\eta(\cdot)$ is the modulus of continuity of the mapping $f$, then $q \leqslant \eta(d)$, where $d$ is the diameter of the covering $C$. There is a natural relation between the admissible paths on the symbolic image $G$ and the trajectories of the dynamical system. The next theorem describes the relationship between $\varepsilon$-trajectories and admissible paths on $G$ and reveals their connection to parameters of the symbolic image.

**Theorem 1** (see [1]). 1. Let $\{z_k\}$ be an admissible path on the symbolic image $G$. Then there exists a sequence of points $\{x_k\}$, $x_k \in M(z_k)$, which is an $\varepsilon$-trajectory for $f$ for any $\varepsilon > d$. In particular, if the sequence $\{z_1, z_2, \ldots, z_p = z_0\}$ is p-periodic, then so is the $\varepsilon$-trajectory $\{x_1, x_2, \ldots, x_p = x_0\}$.

2. Let $\{z_k\}$ be an admissible path on the symbolic image $G$ and let $x_k \in M(z_k)$. Then $\{x_k\}$ is an $\varepsilon$-trajectory for $f$ for any $\varepsilon > q + d$.

3. *There exists $r > 0$ such that if a sequence of points $\{x_k\}$ is an $\varepsilon$-trajectory for $f$, where $\varepsilon < r$, and if $x_k \in M(z_k)$, then the sequence $\{z_k\}$ is an admissible path on the symbolic image $G$. In particular, if the $\varepsilon$-trajectory $\{x_1, x_2, \ldots, x_p = x_0\}$ is $p$-periodic, then $\{z_1, z_2, \ldots, z_p = z_0\}$ is a $p$-periodic path on $G$.*

If $\zeta = \{x_k,\ z \in \mathbb{Z}\}$ is a pseudotrajectory, then the sequence $\omega = \{i_k\colon x_k \in M(i_k)\}$ will be called the trace of the pseudotrajectory $\zeta$ on the symbolic image $G$. Theorem 1 describes the conditions for an admissible path to be the trace of an $\varepsilon$-trajectory, and vice versa.

So, paths on a symbolic image can also be looked upon as encodings of pseudotrajectories.

**Definition 4.** A vertex of a symbolic image is said to be *recurrent* if there is a periodic path through it. Two recurrent vertices $i$ and $j$ are called *equivalent* if there exists a periodic path passing through $i$ and $j$.

We denote the set of recurrent vertices by $RV$. The set of recurrent vertices $RV$ can be partitioned into equivalence classes. In the terminology of graph theory, classes of equivalent recurrent vertices are strongly connected components. In [1] it was shown that the set
$$P(d) = \bigcup_{i \in RV} M(i)$$
is a closed neighbourhood of a chain recurrent set converging to this set as $d \to 0$. There exist depth-first search algorithms for constructing the set $RV$ and the strongly connected components (see [20] and [21]). This makes it possible to localize a chain recurrent set numerically (see [1]).

**The subdivision procedure.** A sequence of symbolic images will be constructed by subdividing coverings. The principal step of the subdivision procedure is as follows. Let $C = \{M(i)\}$ be a covering and let $G$ be the symbolic image for $C$. Assume that a new covering $NC$ is a subdivision of the covering $C$. This means that each cell $M(i)$ is subdivided into cells $m(i, k)$, $k = 1, 2, \ldots$; that is,
$$\bigcup_k m(i, k) = M(i).$$
We let $NG$ denote the new symbolic image for the covering $NC = \{m(i, k)\}$. The vertices of $NG$ will be denoted by $(i, k)$. This construction defines a single-valued mapping $s$ from $NG$ onto $G$, which transforms a vertex $(i, k)$ to the vertex $i$; that is, $s(i, k) = i$. Since
$$f(m(i, k)) \cap m(j, l) \neq \varnothing$$
for small cells, the analogous intersection of the image of the large cel $f(M(i))$ and $M(j)$ is also nonempty:
$$f(M(i)) \cap M(j) \neq \varnothing.$$
Hence $s$ transforms the arc $(i, k) \to (j, l)$ into the arc $i \to j$. Therefore, $s$ maps the directed graph $NG$ onto the directed graph $G$. Since $s$ is a mapping of directed graphs, it takes each admissible path on $NG$ onto some admissible path on $G$.

In particular, the image of a periodic path is a periodic path, and the image of a recurrent vertex is a recurrent vertex. Moreover, the image of the class $NH$ of equivalent recurrent vertices (on $NG$) lies in the class $H$ of equivalent recurrent vertices on $G$.

Many of the properties of a sequence of symbolic images can be conveniently formulated in terms of diagrams. We let $V$ and $NV$ denote the vertex sets of the graphs $G$ and $NG$, respectively. As we pointed out above, the graph $G$ can be interpreted as a set-valued mapping from $V$ onto $V$ such that $G(i) = \{j \colon i \to j\}$. Similarly, the graph $NG$ is a set-valued mapping $NG \colon NV \to NV$. So we have the commutative diagram

$$\begin{array}{ccc} NV & \xrightarrow{NG} & NV \\ \downarrow{s} & & \downarrow{s} \\ V & \xrightarrow{G} & V \end{array} \qquad (2.1)$$

where the commutativity has the same meaning as before:

$$s(NG(i,k)) \subset G(s(i,k)). \qquad (2.2)$$

A recurrent trajectory is defined differently by different authors. For example, Definition 3.3.2 of a recurrent trajectory in [13] is the same as that of a Poisson-stable trajectory in [14]. The paper [22] studies the detailed relationships between these concepts. In this paper we use the definition of a recurrent trajectory in the sense of Birkhoff.

**Definition 5.** A trajectory $K$ is called *recurrent* if for each $\varepsilon > 0$ there exists a positive integer $p$ such that any interval of length $p$ of this trajectory approximates the whole of the trajectory $K$ with accuracy $\varepsilon$.

Consider a sequence of coverings $\{C_t, t \in N\}$ of the manifold $M$ by cells obtained by successive subdivisions; that is, cells in the covering $C_{t+1}$ are obtained by subdividing cells in the covering $C_t$. Let the diameter of the covering $C_t$ tend to zero as $t \to \infty$. Consider the sequence $\{G_t\}$ of symbolic images of the mapping $f \colon M \to M$ with respect to the coverings $C_t$. This gives us the commutative diagram

$$\begin{array}{ccccccc} V_1 & \xleftarrow{s_1} & V_2 & \xleftarrow{s_2} & V_3 & \xleftarrow{s_3} & \cdots \\ \downarrow{G_1} & & \downarrow{G_2} & & \downarrow{G_3} & & \\ V_1 & \xleftarrow{s_1} & V_2 & \xleftarrow{s_2} & V_3 & \xleftarrow{s_3} & \cdots \end{array} \qquad (2.3)$$

In what follows we drop the index $t$ in the notation for $s_t$ when no confusion can ensue.

Consider the limiting case when the cell consists of one point, $M(x) = \{x\}$. In this case the vertex set is the point set $M$ equipped with the discrete topology. The arc set is the family of pairs of the form $(x, f(x))$. We can regard the symbolic image thus defined as coinciding with the original mapping $f \colon M \to M$. Let $G$ be the symbolic image for the covering $C = \{M(i)\}$ and $V$ the vertex set. Above, we defined the set-valued mapping $h \colon M \to V$, which associates with a point $x$ all vertices $i$ such that $x \in M(i)$. The following diagram is commutative:

$$\begin{array}{ccc} V & \xleftarrow{h} & M \\ G \downarrow & & \downarrow f \\ V & \xleftarrow{h} & M \end{array} \qquad (2.4)$$

Taking this into account, we can complete the diagram (2.3) as

$$\begin{array}{ccccccc} V_1 & \xleftarrow{s} & V_2 & \xleftarrow{s} & V_3 & \xleftarrow{s} \cdots \xleftarrow{h} & M \\ G_1 \downarrow & & G_2 \downarrow & & G_3 \downarrow & & \downarrow f \\ V_1 & \xleftarrow{s} & V_2 & \xleftarrow{s} & V_3 & \xleftarrow{s} \cdots \xleftarrow{h} & M \end{array} \qquad (2.5)$$

Each $s$ is a mapping of directed graphs; it maps an admissible path to an admissible path. By the definition of a subdivision, if $s(i) = j$, then the cell $M(i)$ lies in the subdivision of the cell $M(j)$ and $M(i) \subset M(j)$. Let $\omega_t = \{z_k^t,\ k \in \mathbb{Z}\}$ be a path on the symbolic image $G_t$. We denote the space of paths on $G_t$ by $P_t$. The mapping $s\colon G_{t+1} \to G_t$ of directed graphs defines the mapping of paths

$$s\colon P_{t+1} \to P_t, \qquad s(P_{t+1}) \subset P_t,$$

but in general $s(P_{t+1}) \neq P_t$. Fixing a path $\omega_t$ on each symbolic image $G_t$, we get a sequence of paths $\{\omega_t \in P_t\}$, whose elements are in general unrelated. However, by Theorem 1, any trajectory $T = \{x_k = f^k(x_0),\ k \in \mathbb{Z}\}$ defines an admissible path $\omega_t = \{z_k^t \colon x_k \in M(z_k^t)\}$ on each $G_t$ and these paths are related via the mappings $s$:

$$\omega_t = s(\omega_{t+1}). \qquad (2.6)$$

The above path $\omega_t$ can be interpreted as an encoding of the trajectory $T$ with respect to the covering $C_t$.

**Definition 6.** Let $\{C_t\}$ be a sequence of closed coverings, each of which is obtained by subdividing the previous covering. A sequence of admissible paths $\{\omega_t \in P_t\}$ is called *consistent* if (2.6) holds for each $t$.

**Theorem 2.** *Let $\{C_t\}$ be a sequence of closed coverings each of which is obtained by subdividing the previous covering, let their diameters $d_t$ tend to zero, and suppose that on each $G_t$ a path $\omega_t = \{i_k^t,\ k \in \mathbb{Z}\}$ is defined such that the sequence of paths $\{\omega_t\}$ is consistent. Then the following results hold.*

*1. There exists a unique trajectory $T = \{x_k \colon x_{k+1} = f(x_k)\}$ such that $x_k \in M(i_k^t)$ for any $t$.*

*2. If points $x_k^t$ lie in the cells $M(i_k^t)$, then for each $k$ the limit $\lim_{t\to\infty} x_k^t$ exists and coincides with $x_k$.*

*3. If each path $\omega_t$ is periodic, then the shadowed trajectory $T$ is recurrent.*

*Proof.* 1. Let $\omega_t = \{i_k^t,\ k \in \mathbb{Z}\}$. Given a fixed $k$, consider the sequence of cells $\{M(i_k^t),\ t = 1, 2, \dots\}$. Since $s(i_k^{t+1}) = i_k^t$, the cell $M(i_k^{t+1})$ belongs to the subdivision of $M(i_k^t)$. Hence we have

$$M(i_k^1) \supset M(i_k^2) \supset \cdots \supset M(i_k^t) \supset M(i_k^{t+1}) \supset \cdots. \qquad (2.7)$$

Since the cells are closed and their diameters tend to zero as $d_t \to 0$, there exists a unique point
$$x_k = \lim_{t \to \infty} M(i_k^t) = \bigcap_t M(i_k^t).$$
Since the path $\omega_t = \{i_k^t, \ k \in \mathbb{Z}\}$ is admissible, we have
$$f(M(i_k^t)) \cap M(i_{k+1}^t) \neq \varnothing$$
for any $t$. From (2.7) and $d_t \to 0$ it follows that the closed sets $\{f(M(i_k^t)), l = 1, 2, \ldots\}$ are nested and their diameters $q_t$ tend to zero. Hence there exists a unique point
$$x_{k+1}^* = \lim_{t \to \infty} f(M(i_k^t)).$$
A similar analysis shows that the sequence of closed sets $\{f(M(i_k^t)) \cap M(i_{k+1}^t), t = 1, 2, \ldots\}$ has a limit point
$$\lim_{t \to \infty} f(M(i_k^t)) \cap M(i_{k+1}^t) = x_{k+1}^{**}.$$
Since
$$x_{k+1}^{**} \in f(M(i_k^t)) \cap M(i_{k+1}^t) \subset M(i_{k+1}^t)$$
for any $t$, it follows that $x_{k+1}^{**}$ lies in $M(i_{k+1}^t)$, and therefore it lies in
$$\bigcap_t M(i_{k+1}^t) = x_{k+1}.$$
Hence $x_{k+1}^{**} = x_{k+1}$. Since
$$x_{k+1}^{**} \in f(M(i_k^t)) \cap M(i_{k+1}^t) \subset f(M(i_k^t))$$
for any $t$, we see that $x_{k+1}^{**}$ belongs to $f(M(i_k^t))$, and therefore $x_{k+1}^{**}$ lies in
$$\bigcap_t f(M(i_k^t)) = x_{k+1}^*,$$
that is, $x_{k+1}^{**} = x_{k+1}^*$. Next, $f(x_k) \in f(M(i_k^t))$ for any $t$, and hence $f(x_k)$ lies in
$$\bigcap_t f(M(i_k^t)) = x_{k+1}^*,$$
that is, $f(x_k) = x_{k+1}$. This proves assertion 1 of the theorem.

2. In the proof of assertion 1 we constructed the points
$$x_k = \lim_{t \to \infty} M(i_k^t) = \bigcap_t M(i_k^t)$$
and showed that $x_{k+1} = f(x_k)$. Assume that the points $x_k^t$ lie in the cells $M(i_k^t)$. Given $k$, we show that each sequence $\{x_k^t : t \in \mathbb{N}\}$ converges to $x_k$. By assumption $x_k^t \in M(i_k^t)$, and by the above construction $x_k \in M(i_k^t)$. Hence
$$\rho(x_k, x_k^t) \leqslant \operatorname{diam} M(i_k^t) \leqslant d_t \to 0$$
and $\lim_{t \to \infty} x_k^t = x_k$, verifying assertion 2.

3. Assume that each path $\omega_t = \{i_k^t,\ k \in \mathbb{Z}\}$ is periodic with period $p_t$; that is, $i_{k+p_t}^t = i_k^t$ for any $k \in \mathbb{Z}$. Consider the path $\omega_1$ on the symbolic image $G_1$. By the proof of assertion 1 of the theorem, the trajectory $T = \{x_k \colon x_{k+1} = f(x_k)\}$ is contained in the union of the cells

$$U_1 = \left\{ \bigcup_k M(i_k^1) \colon i_k^1 \in \omega_1 \right\}.$$

By hypothesis $\omega_1$ is a periodic path with period $p_1$. Hence $U_1$ is a union of a finite number of cells $M(i_k^1)$, $k = 1, 2, \ldots, p_1$. Let $d_1$ be the diameter of the covering $C_1$. Then the ball $B_k^1$ with radius $d_1$ and centre $x_k \in M(i_k^1)$ contains the cell $M(i_k^1)$. Therefore, the union of the balls

$$\bigcup B_k^1, \qquad k = 1, 2, \ldots, p_1,$$

contains the trajectory $T$. Hence the $d_1$-neighbourhood of any interval of $T$ of length $p_1$ contains the whole of $T$. Let $d_m$ be the diameter of the covering $C_m$ and let $\omega_m$ be the path of period $p_m$ on the symbolic image $G_m$. A similar argument for the covering $C_m$ shows that the union of the balls $B_k^m$ of radius $d_m$ with centres at the points $x_k \in M(i_k^m)$, $k = 1, 2, \ldots, p_m$, contains the trajectory $T$. Therefore, the $d_m$-neighbourhood of any interval of $T$ of length $p_m$ contains the whole of $T$. The trajectory $T$ is recurrent because, for any fixed positive $\varepsilon$, we can find $d_m < \varepsilon$ and a period $p_m$ such that the $d_m$-neighbourhood of any interval of $T$ of length $p_m$ contains the whole of $T$. This proves assertion 3, and thus also Theorem 2.

*Remark.* By Theorem 2, for the sequence of paths $\{\omega_l\}$ described in the hypotheses of the theorem, there exists a unique trajectory $T$. The converse is not true; that is, there exists a trajectory generating more than one sequence $\{\omega_l\}$ of this form. For example, a trajectory $T$ passing through a boundary point of cells generates more than one sequence of admissible encodings.

Consider a directed graph $G$ with vertices $V = \{i, j, \ldots\}$. Let $\omega = (i_k,\ k \in \mathbb{Z})$ be a path on $G$ and let $P = \{\omega\}$ be the space of paths equipped with the Tychonoff topology. We can introduce the distance on $P$ as follows, for example:

$$\rho(\omega_1, \omega_2) = \sum_{-\infty}^{\infty} \lambda^{|r|} \delta(i_r, j_r);$$

here $\omega_1 = (i_r,\ r \in \mathbb{Z})$, $\omega_2 = (j_r,\ r \in \mathbb{Z})$; $\delta(i, j) = 1$ for $i \neq j$; otherwise $\delta(i, j) = 0$; $0 < \lambda < 1$. The space $P$ is compact in this topology. This topology is also generated by the (different) distance

$$\rho_0(\omega_1,\ \omega_2) = \begin{cases} 2^{-k} & \text{if } k = \max\{m \colon i_r = j_r \text{ for } |r| \leqslant m\}, \\ 0 & \text{if } \omega_1 = \omega_2. \end{cases}$$

This means that the sequences $\omega_1$ and $\omega_2$ coincide for $|r| \leqslant m$, but $i_{m+1} \neq j_{m+1}$ or $i_{-m-1} \neq j_{-m-1}$. For the topological properties of abstract sequence spaces equipped with the Tychonoff topology, see [4].

Let $\{C_t,\ t \in \mathbb{N}\}$ be a sequence of coverings of the manifold $M$ by cells obtained by successive subdivisions and let their diameters $d_t$ tend to zero. We denote the space of paths on $G_t$ equipped with the Tychonoff topology by $P_t$. The mapping $s\colon G_{t+1} \to G_t$ of directed graphs defines a projection in the space of paths:

$$s\colon P_{t+1} \to P_t.$$

This mapping is continuous in the Tychonoff topology, and its image $s(P_{t+1})$ is a compact set in $P_t$. Consider two subdivisions $C_{t+m}$, $m > 0$, and $C_t$. By construction, $C_{t+m}$ is a subdivision of the covering $C_t$. This gives us the mapping $s\colon G_{t+m} \to G_t$ of directed graphs and the mapping $s\colon P_{t+m} \to P_t$ of the space of paths. Assume that on each symbolic image $G_t$ we are given a path $\omega_t$. Then on each $G_t$ we obtain a sequence of paths $\omega_t^m = s(\omega_{t+m})$, $m = 1, 2, \ldots$; that is, $\omega_t^m$ is the projection of the path $\omega_{t+m}$ via the mapping $s\colon P_{t+m} \to P_t$.

**Definition 7.** We say that a sequence $\omega_t$ is *weakly convergent* if the sequence $\omega_t^m$ converges in the Tychonoff topology as $m \to \infty$ for each fixed $t$.

**Theorem 3.** *Let $\{C_t,\ t \in \mathbb{N}\}$ be a sequence of coverings of the manifold $M$ by cells obtained by successive subdivisions, let their diameters $d_t$ tend to zero, and suppose that on each $G_t$ a path $\omega_t = \{i_k^t,\ k \in \mathbb{Z}\}$ is given. If the sequence $\omega_t$ is weakly convergent and if $x_k^t \in M(i_k^t)$, then the limit $\lim_{t\to\infty} x_k^t = x_k$ exists and the sequence $\{x_k,\ k \in \mathbb{Z}\}$ is a trajectory of the mapping $f$.*

*Proof.* The space $P_t$ is compact in the Tychonoff topology. Each path $\omega_t^m \subset P_t$ is the projection of $\omega_{t+m}$ by means of the mapping $s\colon P_{t+m} \to P_t$. By the hypothesis of the theorem, the sequence $\{\omega_t^m,\ m \in \mathbb{N}\}$ converges as $m \to \infty$. We set

$$\lim_{m\to\infty} \omega_t^m = \omega_t^*.$$

We will show that

$$s(\omega_{t+1}^*) = \omega_t^*.$$

In fact, $\omega_t^m = s(\omega_{t+m})$ and $\omega_{t+1}^m = s(\omega_{t+1+m})$. Since $s$ is continuous (in the Tychonoff topology), we have

$$\lim_{m\to\infty} \omega_t^m = \lim_{m\to\infty} s(\omega_{t+m}) = \lim_{r\to\infty} s(\omega_{t+1+r}) = \lim_{r\to\infty} s(\omega_{t+1}^r) = s\left(\lim_{r\to\infty} \omega_{t+1}^r\right).$$

This implies the required equality; that is, the sequence $\omega_t^*$ is consistent. Let $\omega_t^*$ be the sequence $\{i_k^{*t},\ k \in \mathbb{Z}\}$ and let $x_k^{*t}$ be a point in the cell $M(i_k^{*t})$. By assertion 2 of Theorem 2, for each $k$, the limit

$$\lim_{t\to\infty} x_k^{*t} = x_k$$

exists and the sequence $\{x_k,\ k \in \mathbb{Z}\}$ is a trajectory $T = \{x_k\colon x_{k+1} = f(x_k)\}$.

We show that $\lim_{t\to\infty} x_k^t = x_k$. By the hypothesis of the theorem, each point $x_k^t$ lies in the cell $M(i_k^t)$, where $\omega_t = \{i_k^t,\ k \in \mathbb{Z}\}$. We fix $k$ and $\varepsilon > 0$. Given $\varepsilon$, we find $\tau$ for which the diameter of the covering $d_\tau$ is smaller than $\varepsilon$. For this $\tau$ we find a number $m(l)$ such that the sequences $\omega_\tau^*$ and $\omega_\tau^m$ are sufficiently close in the

Tychonoff topology for their elements with indices $k$ such that $|k| \leqslant l$ to be equal; that is,
$$i_k^{*\tau} = s(i_k^{\tau+m}), \qquad |k| \leqslant l,$$
for all $m > m(l)$. The equality
$$s(\omega_{t+1}^*) = \omega_t^*$$
implies that
$$s(\omega_{t+m}^*) = \omega_t^*.$$
Hence, for the projections of the indices of the cells, we have
$$s(i_k^{*\tau+m}) = i_k^{*\tau} = s(i_k^{\tau+m}),$$
where $s$ maps the directed graph $G_{\tau+m}$ onto $G_\tau$. This means that the points $x_k^{*\tau+m}$, $x_k^{\tau+m}$ lie in the same cell $M(i_k^{*\tau})$, which has diameter strictly less than $\varepsilon$. Hence, for $t > \tau + m$,
$$\rho(x_k^t, x_k^{*t}) < \varepsilon, \qquad |k| < l.$$
Since $\varepsilon$ is an arbitrary positive number, the equality
$$\lim_{t \to \infty} x_k^{*t} = x_k$$
implies that
$$\lim_{t \to \infty} x_k^t = x_k.$$
Theorem 3 is proved.

**The sequence space on a manifold.** Let $\Lambda$ be the set of doubly infinite sequences $\zeta = \{x_k, \ k \in \mathbb{Z}\}$ of points of the manifold $M$. Here we view $x_k$ as the $k$th coordinate of the sequence $\zeta$. We define a distance $r$ on $\Lambda$ as follows. Let $\zeta_1 = \{x_k, \ k \in \mathbb{Z}\}$ and $\zeta_2 = \{y_k, \ k \in \mathbb{Z}\}$, and let $\rho(x, y)$ be the distance on $M$. Then we set
$$r(\zeta_1, \zeta_2) = \sum_{-\infty}^{\infty} \frac{\rho(x_k, y_k)}{2^{|k|}}.$$

The distance $r$ is an analogue of the distance between paths on the graph, which generates the Tychonoff topology in the space of paths. The above topology on the space of sequences $\Lambda$ will be called the *Tychonoff topology*.

In the proof of assertion 2 of Theorem 2, the sequence $\{x_k^t\}$ converges to the sequence $\{x_k\}$ uniformly in $k$, because $\rho(x_k^t, x_k) \leqslant d_t$. In the proof of Theorem 3 the sequence $\{x_k^t\}$ converges to the sequence $\{x_k\}$ in the Tychonoff topology of the space of sequences $\Lambda$, because the paths $\omega_\tau^*$ and $\omega_\tau^m$ are close in the Tychonoff topology.

By the definition of the distance $r$,
$$\rho(x_k, y_k) < 2^{|k|} r(\zeta_1, \zeta_2).$$

Therefore, if $\zeta_1$ converges to $\zeta_2$ in the Tychonoff topology of the space $\Lambda$, then $x_k$ converges to $y_k$ on $M$; that is, the sequence $\zeta_1$ converges to $\zeta_2$ coordinatewise. The converse result also holds.

**Proposition 1.** *Given a compact manifold $M$, if sequences converge coordinate-wise then they converge in the Tychonoff topology.*

**Corollary 1.** *The set of sequences $\Lambda$ is compact in the Tychonoff topology.*

*Proof.* Let $\zeta_t = \{x_k^t, \ k \in \mathbb{Z}\}$ be a sequence of elements of $\Lambda$. Since $M$ is compact, for $k = 0$ there exists a subsequence $\tau = t_p \to \infty$ with convergent 0th coordinates: $x_0^\tau \to x_0^*$. For $k = \pm 1$, we select a subsequence from the sequence $\zeta_\tau$ with convergent 1st and $(-1)$th coordinates: $x_{\pm 1}^\tau \to x_{\pm 1}^*$. In a similar way, we construct a subsequence with convergent $k$th and $(-k)$th coordinates: $x_{\pm k}^\tau \to x_{\pm k}^*$. We then construct a diagonal sequence which converges coordinatewise.

This proves Corollary 1.

**Proposition 2.** *If $\zeta_t$ is a sequence of $\varepsilon_t$-trajectories and $\varepsilon_t$ tends to zero, then there exists a subsequence $\zeta_{t_m}$ converging to a trajectory in the Tychonoff topology.*

*Proof.* Let $x_0^t$ be the zeroth element of the sequence $\zeta_t$. The points $\{x_0^t, \ t \in \mathbb{N}\}$ lie in a compact manifold, and hence there exists a subsequence $x_0^\tau$ indexed by $\tau = t_m$ which converges to $x^*$. We show that the trajectory through the point $x^*$ is the one we require.

Consider the $\varepsilon$-trajectory $\zeta = \{x_k, \ k \in \mathbb{Z}\}$. We estimate the distance $\rho(x_k, f^k(x_0))$ in terms of $\varepsilon$ and the modulus of continuity $\eta$ of the mapping $f$: $\rho(f(x), f(y)) \leqslant \eta(\rho(x, y))$, where $\eta$ can be assumed to be monotonic. For $k > 0$ we have

$$\rho(x_k, f^k(x_0)) \leqslant \rho(x_k, f(x_{k-1})) + \rho(f(x_{k-1}), f^2(x_{k-2})) \\ + \rho(f^2(x_{k-2}), f^3(x_{k-3})) + \cdots + \rho(f^{k-1}(x_1), f^k(x_0)).$$

We can estimate each term:

$$\rho(x_k, f(x_{k-1})) \leqslant \varepsilon,$$
$$\rho(f(x_{k-1}), f^2(x_{k-2})) \leqslant \eta(\rho(x_{k-1}, f(x_{k-2}))) \leqslant \eta(\varepsilon),$$
$$\rho(f^2(x_{k-2}), f^3(x_{k-3})) \leqslant \eta(\rho(f(x_{k-2}), f^2(x_{k-3}))) \leqslant \eta^2(\varepsilon),$$

where $\eta^2$ is the composition of the function $\eta$ with itself. A similar analysis shows that

$$\rho(f^{k-1}(x_1), f^k(x_0)) \leqslant \eta^{k-1}(\varepsilon).$$

So,

$$\rho(x_k, f^k(x_0)) \leqslant \sum_{l=0}^{k-1} \eta^l(\varepsilon).$$

Hence, for fixed $k$,

$$\rho(x_k, f^k(x_0)) \to 0$$

as $\varepsilon \to 0$. Using the modulus of continuity of $f^{-1}$ we can show that the last estimate holds for $k \leqslant 0$. By Proposition 1, the distance in the Tychonoff topology tends to zero,

$$r(\zeta_\tau, T(x_0^\tau)) \to 0.$$

Since the homeomorphism $f$ is continuous, the trajectory $T(x_0^\tau)$ converges to the trajectory $T(x^*)$ in the Tychonoff topology, because the sequence $x_0^\tau$ converges to $x^*$. Hence $\zeta_\tau$ converges to $T(x^*)$ in the Tychonoff topology.

Proposition 2 is proved.

**Corollary 2.** *The set of trajectories* $\mathrm{Tr} \subset \Lambda$ *is a compact set in the Tychonoff topology.*

**Theorem 4.** *Let* $\{C_t,\ t \in \mathbb{N}\}$ *be a sequence of coverings of the manifold $M$ by cells obtained by successive subdivisions with diameters* $d_t \to 0$ *and suppose that a path* $\omega_t = \{i_k^t,\ k \in \mathbb{Z}\}$ *is defined on each $G_t$. Then there exists a subsequence indexed by* $\tau = t_m \to \infty$ *such that, for any* $x_k^\tau \in M(i_k^\tau)$, *the limit*

$$\lim_{\tau \to \infty} x_k^\tau = x_k$$

*exists and the sequence* $\{x_k,\ k \in \mathbb{Z}\}$ *is a trajectory of the mapping $f$.*

*Proof.* Consider the symbolic image $G_1$ and the space of paths $P_1$. We project all the paths $\omega_t$, $t > 1$, onto $P_1$ by mappings of the form $s$. This produces the sequence of paths $\omega_t^1$ in $P_1$. Since $P_1$ is a compact set, there exists a convergent subsequence $\omega_m^1$, $m = t_r \to \infty$. Next, we consider the symbolic image $G_2$ and the space of paths $P_2$, and project all the selected paths $\omega_m$ onto $P_2$. As in the previous construction, we construct a convergent subsequence $\omega_l^2$. Proceeding in this way, on each $G_t$ we construct a sequence $\omega_*^t$, which is a subsequence of the previous sequence in the following sense: the projection $s\omega_*^t$ is a subsequence of $\omega_*^{t-1}$. If we take the diagonal sequence $\omega_\tau^\tau$, we obtain a sequence of paths convergent in the Tychonoff topology. Now the conclusion of the theorem follows from Theorem 3.

In all the above results we assumed that the diameter of the coverings tends to zero. However, for practical purposes, it is desirable to have a result describing the shadowing by paths of a symbolic image for a covering with sufficiently small positive diameter. The next theorem, Theorem 5, gives the basis for a practical application of symbolic images. Recall that if $\omega = \{i_k,\ k \in \mathbb{Z}\}$ is an admissible path on a symbolic image $G$, then the sequence of points $\{x_k,\ x_k \in M(i_k)\}$ is a pseudotrajectory (by Theorem 1), which is the trace of the path $\omega$.

**Theorem 5.** *Let $U$ be a neighbourhood of the set* $\mathrm{Tr}$ *in the Tychonoff topology. Then there exists* $d_0 > 0$ *such that, for any symbolic image $G$, constructed for the covering $C$ with diameter* $d \leqslant d_0$, *any trace $\chi(\omega)$ of any admissible path $\omega$ for $G$ lies in $U$.*

*Proof.* Assume by contradiction that there exists a neighbourhood $U$ of the set of trajectories $\mathrm{Tr}$ for which the conclusion of the theorem does not hold. Then, for any $d > 0$ there is a covering $C$ of diameter $d$ such that the symbolic image $G(C)$ contains an admissible path $\omega = \{i_k\}$ for which the sequence $\{x_k,\ x_k \in M(i_k)\}$ lies outside $U$. According to Theorem 1, the sequence $\{x_k,\ x_k \in M(i_k)\}$ is an $\varepsilon$-trajectory of the mapping $f$ for any $\varepsilon > q + d$. We can take $\varepsilon = \eta(d) + d$, where $\eta(\cdot)$ is the modulus of continuity of $f$. Let $C_t$ be a sequence of coverings such that the subdivision diameters $d_t$ tend to 0 and the above requirements are satisfied. As a result, we obtain the sequence $\zeta_t = \{x_k^t\}$ of $\varepsilon_t$-trajectories with $\varepsilon_t \to 0$ that lie outside $U$. By Proposition 2, there exists a subsequence $\zeta_{t_m}$ converging to the trajectory $T$ in the Tychonoff topology. But then $T$ does not lie in Tr. This contradiction completes the proof.

Consider the symbolic image $G$ of the mapping $f$. The set of (admissible) paths $P$ on $G$ is equipped with the Tychonoff topology; Theorem 1 describes the relationship between the paths on $G$ and the pseudotrajectories of the mapping $f$. Any trajectory $T = \{x_k, \ k \in \mathbb{Z}\}$ has its trace on $G$ in the form of the path $\omega = \{i_k \colon x_k \in M(i_k)\}$, which can naturally be called an encoding of the trajectory $T$. The set-valued mapping $h(x) = \{i \colon x \in M(i)\}$ defines an encoding of any trajectory. We let Cod denote the set of paths encoding trajectories. Our next aim is to describe the embedding $\mathrm{Cod} \subset P$.

Let $\{C_t, \ t \in \mathbb{N}\}$ be a sequence of coverings of the manifold $M$ by cells obtained by successive subdivisions and let their diameters $d_t$ tend to zero. These coverings generate a sequence of symbolic images $G_t$ and a set of paths $P_t$. The set of encodings $\mathrm{Cod}_t$ lies in $P_t$. The subdivision of cells defines a mapping $s \colon G_{t+1} \to G_t$ of directed graphs and a mapping of paths $s \colon P_{t+1} \to P_t$. The projection $s(P_{t+1})$ is a compact set in $P_t$. Since the covering $C_{t+m}$ is a subdivision of $C_t$ for any $m > 0$, there exists a mapping $s \colon G_{t+m} \to G_t$ of directed graphs and a projection $s \colon P_{t+m} \to P_t$ of paths. This gives the sequence of projections $P_t^m = s(P_{t+m})$ for $m > 0$ in the space of paths $P_t$.

**Theorem 6.** *Let $\{C_t, \ t \in \mathbb{N}\}$ be a sequence of coverings of the manifold $M$ by cells obtained by successive subdivisions, and let their diameters $d_t$ tend to zero. Then:*

- *the projections $\{P_t^m, \ m \in \mathbb{N}\}$ form a sequence of nested compact sets; that is,*

$$P_t \supset P_t^1 \supset P_t^2 \supset \cdots \supset P_t^m \supset P_t^{m+1} \cdots ;$$

- *the set of encodings $\mathrm{Cod}_t$ is a compact subset of $P_t$ and*

$$\mathrm{Cod}_t = \lim_{m \to \infty} P_t^m = \bigcap_{m > 0} P_t^m;$$

- *the $s$-image of the set of encodings of trajectories coincides with the encodings of trajectories; that is,*

$$s(\mathrm{Cod}_{t+1}) = \mathrm{Cod}_t.$$

*Proof.* 1. The mapping $s \colon G_{t+1} \to G_t$ of directed graphs defines a mapping of paths $s \colon P_{t+1} \to P_t$, which is continuous in the Tychonoff topology. The image $s(P_{t+1}) = P_t^1$ is a compact set in $P_t$. The mapping $s \colon P_{t+m} \to P_t$ is a composition of the above mappings $s \colon P_{\tau+1} \to P_\tau$. Hence all the projections $P_t^m = s(P_{t+m})$, $m \in \mathbb{N}$, are compact. Since $s(P_{t+1}) \subset P_t$ for any $t$, $\{P_t^m, \ m \in \mathbb{N}\}$ is a sequence of nested compact sets.

2. For each $t$ we have $\mathrm{Cod}_t \subset P_t$ and $s(P_{t+1}) \subset P_t$. Hence $\mathrm{Cod}_t \subset P_t^m$ for any $m > 0$, and therefore,

$$\mathrm{Cod}_t \subset \bigcap_{m>0} P_t^m.$$

We prove the reverse inclusion. Let $\omega \in \bigcap_{m>0} P_t^m$; then there exists a sequence $\omega_m \in P_{t+m}$ such that $s(\omega_m) = \omega$. Therefore, the sequence $\omega_m$ is consistent, and so, by Theorem 2, there exists a trajectory $T = \{x_k\}$ which is the trace of each path $\omega_m$ and of $\omega$. In other words, $\omega_m$ and $\omega$ are encodings of the trajectory $T = \{x_k\}$, which proves the inclusion $\omega \in \mathrm{Cod}_t$.

3. The image is found as follows:
$$s(\mathrm{Cod}_{t+1}) = s\left(\bigcap_{m>0} P^m_{t+1}\right) = \bigcap_{m>0} s(P_{t+m}) = \bigcap_{m>0} P^m_t = \mathrm{Cod}_t.$$

The theorem is proved.

## § 3. Shadowing of ergodic measures

We recall that $\mathrm{ext}(A)$ denotes the set of extreme points of a convex set $A$. The set of invariant measures $\mathcal{M}(f)$ is a convex weakly compact set, whose extreme points are ergodic measures. The convex hull of a set $A$ is the intersection of all convex sets that contain $A$. According to Theorem 1.2.2 in [23], the convex hull of $A$ is the set of convex combinations of finite numbers of points in $A$. The convex hull thus defined need not be closed. By $\mathcal{M}^0$ we denote the set consisting of measures each of which is the convex hull of a finite number of ergodic measures,

$$\mathcal{M}^0 = \left\{ \mu = \sum_{k=1}^m \alpha_k \mu_k,\ \mu_k \in \mathrm{ext}\,\mathcal{M}(f),\ \alpha_k \geqslant 0,\ \sum_k \alpha_k = 1 \right\}.$$

The set $\mathcal{M}^0$ is convex, and moreover,
$$\mathcal{M}^0 \subset \mathcal{M}(f).$$

According to [14], the weak closure of $\mathcal{M}^0$ coincides with $\mathcal{M}(f)$,
$$\overline{\mathcal{M}^0} = \mathcal{M}(f).$$

By Choquet's theorem, any point in a convex compact set is an integral sum of its extreme points. We also need the following result on expanding a system with finite invariant measure in ergodic components (see Theorem 4.1.12 in [13]).

**Theorem 7.** *For each invariant measure $\mu \in \mathcal{M}(f)$, there exists a division (modulo nullsets) of the manifold $M$ into invariant measurable sets $\Omega_\alpha$, $\alpha \in A$, where $A$ is an abstract measure space and each $\Omega_\alpha$ is equipped with an ergodic invariant measure $\mu_\alpha$, such that, for any integrable function $\varphi$,*

$$\int_M \varphi\,d\mu = \int_A \int_M \varphi\,d\mu_\alpha\,d\alpha. \tag{3.1}$$

Let $C$ be a covering of $M$ by polyhedral cells intersecting in boundary discs, and let $C^*$ be the division obtained from $C$ by attributing each boundary disc to one of the adjoining polyhedra. Let $G$ be the symbolic image constructed for the covering $C$. The set-value mapping $h\colon M \to G$, $h(x) = \{i\colon x \in M(i)\}$, is single-valued on $C^*$. This mapping generates the mapping

$$h^*(\mu) = m = \{m_{ij} = \mu(f(M^*(i)) \cap M^*(j))\}$$

of the set of invariant measures $\mathscr{M}(f)$ into the set of flows $\mathscr{M}(G)$, where $\mu \in \mathscr{M}(f)$, $M^*(i) \in C^*$, and $m = \{m_{ij}\}$ is a flow on $G$. The mapping $h^* \colon \mathscr{M}(f) \to \mathscr{M}(G)$ is linear in the following sense:

$$h^*(\alpha \mu_1 + \beta \mu_2) = \alpha h^*(\mu_1) + \beta h^*(\mu_2).$$

Note that $h^* \colon \mathscr{M}(f) \to \mathscr{M}(G)$ is not continuous as a mapping of the space $\mathscr{M}(f)$ equipped with the weak topology into the space $\mathscr{M}(G)$ with the metric $\rho(m^1, m^2) = \sum_{ij} |m^1_{ij} - m^2_{ij}|$. By Theorem 2.1 in [24], $h^*$ is continuous at a point $\mu$ if the boundaries of the cells have zero measure.

By definition of an ergodic measure, if an ergodic measure $\mu$ lies in a convex set $W \subset \mathscr{M}(f)$, then $\mu$ is an extreme point of $W$. Recall that the simple flow $m$ is generated by a simple cycle and cannot be decomposed in terms of other flows. If a simple flow $m$ is contained in a convex set $V \subset \mathscr{M}(G)$, then $m$ is an extreme point of $V$.

**Theorem 8.** *Let $h^* \colon \mathscr{M}(f) \to \mathscr{M}(G)$ be a mapping of the set of invariant measures into the set of flows on the graph $G$. Then the images $h^*(\mathscr{M}(f))$ and $h^*(\mathscr{M}^0)$ are convex subsets of $\mathscr{M}(G)$, and moreover, each extreme point of the image $h^*(\mathscr{M}^0)$ is the image of an ergodic measure,*

$$\operatorname{ext}(h^*(\mathscr{M}^0)) \subset h^*(\operatorname{ext}(\mathscr{M}(f))).$$

*Proof.* The mapping $h^* \colon \mathscr{M}(f) \to \mathscr{M}(G)$ is linear, and therefore the images $h^*(\mathscr{M}(f))$ and $h^*(\mathscr{M}^0)$ are convex. The convex sets $\mathscr{M}(f)$ and $\mathscr{M}^0$ have the same extreme points, which are ergodic measures. Each invariant measure $\mu \in \mathscr{M}^0$ expands in a finite sum of ergodic measures:

$$\mu = \sum_z \alpha_z \mu_z.$$

Hence the image of the measure $\mu$ expands in a finite sum of images of ergodic measures:

$$h^* \mu = \sum_z \alpha_z h^* \mu_z.$$

Let $h^* \mu = m$ be an extreme point of the convex set $h^*(\mathscr{M}^0)$. The flow $m$ can be decomposed in images of ergodic measures $\{h^* \mu_z\}$,

$$m = \sum_z \alpha_z h^* \mu_z.$$

For the extreme point $m$ such an expansion is only possible under the condition that for $\alpha_z > 0$ the image $h^* \mu_z$ coincides with $m$; that is, $h^* \mu_z = m$, and the sum of all such $\alpha_z > 0$ is 1. It follows that each extreme point of the image $h^*(\mathscr{M}^0)$ is the image of some extreme point of $\mathscr{M}^0$. Therefore,

$$\operatorname{ext}(h^*(\mathscr{M}^0)) \subset h^*(\operatorname{ext}(\mathscr{M}^0)) = h^*(\operatorname{ext}(\mathscr{M}(f))).$$

Theorem 8 is proved.

**Proposition 3.** Let $h^*\colon \mathscr{M}(f) \to \mathscr{M}(G)$ be a mapping of the set of invariant measures into the set of flows on the graph $G$ and let $m \in \mathscr{M}(G)$. Then:

1) the preimage $(h^*)^{-1}(m) = \mathscr{M}(m)$ is a convex set;
2) $\mathscr{M}^0(m) = \mathscr{M}^0 \cap \mathscr{M}(m)$ is a convex set;
3) if $m$ is an extreme point of the convex polyhedron $\mathscr{M}(G)$ and $\mathscr{M}(m) \neq \varnothing$, then the extreme points of the set $\mathscr{M}^0(m)$ are ergodic measures; that is,

$$\operatorname{ext}(\mathscr{M}^0(m)) \subset \operatorname{ext}(\mathscr{M}(f)).$$

*Proof.* 1) The mapping $h^*\colon \mathscr{M}(f) \to \mathscr{M}(G)$ is linear. If $\mu_1$ and $\mu_2$ lie in $\mathscr{M}(f)$, $h^*(\mu_1) = m = h^*(\mu_2)$ and $\alpha + \beta = 1$, then $h^*(\alpha\mu_1 + \beta\mu_2) = m$; that is, $\mathscr{M}(m)$ is a convex set.

2) The set $\mathscr{M}^0(m)$ is convex as the intersection of the convex sets $\mathscr{M}^0$ and $\mathscr{M}(m)$. If $\mu$ lies in $\mathscr{M}^0(m)$, then $\mu$ can be written as a finite sum of ergodic measures, $\mu = \sum \alpha_k \mu_k$, and $h^*(\mu) = m$.

3) Let the measure $\mu$ lie in $\mathscr{M}^0(m)$. Then $\mu$ can be written as a finite sum of ergodic measures,

$$\mu = \sum \alpha_k \mu_k, \qquad \mu_k \in \operatorname{ext} \mathscr{M}^0 = \operatorname{ext} \mathscr{M}(f).$$

Let $m$ be an extreme point of the polyhedron $\mathscr{M}(Q)$ and let $h^*(\mu) = m$. We claim that all the $\mu_k$ lie in $\mathscr{M}^0(m)$. In fact,

$$h^*\mu = \sum \alpha_k h^*\mu_k, \quad \text{or} \quad m = \sum \alpha_k m_k,$$

where $h^*\mu_k = m_k$. Since $m$ is an extreme point, the equality $m = \sum \alpha_k m_k$ is possible if $m_k = m$ and $\sum_k \alpha_k = 1$ for all $k$ such that $\alpha_k > 0$. This means that $h^*\mu_k = m$ and all ergodic measures $\mu_k$ lie in $\mathscr{M}^0(m) \subset \mathscr{M}(m)$. Let $\mu$ be an extreme point of $\mathscr{M}^0(m)$. Arguing by contradiction, we show that $\mu$ is an ergodic measure. If $\mu$ were not ergodic, then it could be written as a sum of ergodic measures,

$$\mu = \sum \alpha_k \mu_k, \qquad \mu_k \in \operatorname{ext} \mathscr{M}(f).$$

By the above, all the $\mu_k$ lie in $\mathscr{M}^0(m)$, and therefore $\mu$ is not an extreme point of $\mathscr{M}^0(m)$. This contradiction completes the proof.

**Corollary 3.** *If $m \in h^*(\mathscr{M}(f))$ is a simple flow, then $m$ is the image of some ergodic measure $\mu$,*

$$m = h^*(\mu).$$

**Proposition 4** (see [19]). *Let $Q$ and $G$ be directed graphs, let $s\colon Q \to G$ be a mapping of directed graphs and let there exist a flow $m$ on $Q$. Then a flow $m^* = s^*m$ is induced on $G$ such that the measure of the arc $i \to j \in G$ is evaluated as*

$$m^*_{ij} = \sum_{s(p \to q) = i \to j} m_{pq},$$

*where the sum is taken over all arcs $p \to q$ mapped to $i \to j$. If an arc $i \to j$ has no preimage, then $m^*_{ij} = 0$.*

**Corollary 4.** *The mapping $s^*\colon \mathscr{M}(Q) \to \mathscr{M}(G)$ is linear in the following sense: if $m_1$ and $m_2$ lie in $\mathscr{M}(Q)$, $\alpha$ and $\beta$ are nonnegative and $\alpha + \beta = 1$, then $s^*(\alpha m_1 + \beta m_2) = \alpha s^*(m_1) + \beta s^*(m_2)$.*

**Proposition 5.** *Let $\mathscr{M}(Q)$ be the set of flows on a directed graph $Q$, let $G$ be a different directed graph, and let $s\colon Q \to G$ be a mapping between directed graphs. Then the image $s^*(\mathscr{M}(Q))$ is a convex polyhedron in $\mathscr{M}(G)$, and the vertices of the polyhedron $s^*(\mathscr{M}(Q))$ lie in the image of the vertices of the polyhedron $\mathscr{M}(Q)$; that is,*

$$\mathrm{ext}(s^*(\mathscr{M}(Q))) \subset s^*(\mathrm{ext}(\mathscr{M}(Q))).$$

*Proof.* The mapping $s^*\colon \mathscr{M}(Q) \to \mathscr{M}(G)$ is linear and the image $s^*(\mathscr{M}(Q))$ is a convex polyhedron. Simple flows are extreme points of the set $\mathscr{M}(Q)$ and the number of simple flows is finite. In [19] it was shown that each flow $m \in \mathscr{M}(Q)$ can be written as a sum of simple flows,

$$m = \sum_z \alpha_z m_z,$$

where $\mathrm{ext}(\mathscr{M}(Q)) = \{m_z\}$ is the complete system of simple flows in $\mathscr{M}(Q)$. Then the image of the flow $m$ has the form

$$m^* = s^*m = \sum_z \alpha_z s^* m_z, \tag{3.2}$$

that is, any flow on the image $s^*(\mathscr{M}(Q))$ can be decomposed in terms of $\{s^* m_z\}$. Let $m^*$ be an extreme point of the convex set $s^*(\mathscr{M}(Q))$; that is, $m^* \in \mathrm{ext}(s^*(\mathscr{M}(Q)))$. Then $m^*$ expands in a sum of flows of the form $\{s^* m_z\}$,

$$m^* = \sum_z \alpha_z s^* m_z.$$

Since $m^*$ is an extreme point, the expansion (3.2) is possible only if all the images $s^* m_z$ for $\alpha_z > 0$ coincide with $m^*$; that is, $s^* m_z = m^*$, and the sum of all such $\alpha_z > 0$ is 1. Hence each extreme point of the image $s^*(\mathscr{M}(Q))$ is the image of an extreme point of the set $\mathscr{M}(Q)$. Therefore, we have the inclusion

$$\mathrm{ext}(s^*(\mathscr{M}(Q))) \subset s^*(\mathrm{ext}(\mathscr{M}(Q))).$$

Proposition 5 is proved.

Let $s\colon Q \to G$ be a mapping of directed graphs. Consider a flow $m^* \in \mathscr{M}(G)$. There are two cases to consider:

1) $m^*$ is not the image of some flow under the mapping $s^*$; that is,

$$m^* \notin s^*(\mathscr{M}(Q));$$

2) $m^*$ is the image of some flow under the mapping $s^*$; that is,

$$m^* \in s^*(\mathscr{M}(Q)).$$

**Proposition 6.** *Let* $s\colon Q \to G$ *be a mapping of directed graphs, let* $m^* \in \mathcal{M}(G)$, *and let* $\mathcal{M}(m^*) = \{m \in \mathcal{M}(Q)\colon s^*(m) = m^*\}$. *Then:*
- *the preimage* $(s^*)^{-1}(m^*) = \mathcal{M}(m^*)$ *is a convex subset of* $\mathcal{M}(Q)$;
- *if* $m^*$ *is a simple flow on* $G$ *and* $\mathcal{M}(m^*) \neq \varnothing$, *then extreme points of the set* $\mathcal{M}(m^*)$ *are simple flows; that is,*

$$\operatorname{ext}(\mathcal{M}(m^*)) \subset \operatorname{ext}(\mathcal{M}(Q)).$$

*Proof.* 1. The mapping $s^*\colon \mathcal{M}(Q) \to \mathcal{M}(G)$ is linear. If $m_1$ and $m_2$ lie in $\mathcal{M}(Q)$, $s^*(m_1) = m^* = s^*(m_2)$ and $\alpha + \beta = 1$, then $s^*(\alpha m_1 + \beta m_2) = m^*$, and therefore $\mathcal{M}(m^*)$ is a convex set.

2. Let $m$ be an extreme point of the set $\mathcal{M}(m^*) \subset \mathcal{M}(Q)$. The point $m$ can be written as a sum of simple flows of the set $\mathcal{M}(Q)$,

$$m = \sum \alpha_k m_k, \qquad m_k \in \operatorname{ext}\mathcal{M}(Q).$$

Hence

$$s^*m = \sum \alpha_k s^* m_k, \quad \text{or} \quad m^* = \sum \alpha_k m_k^*,$$

where $s^* m_k = m_k^*$. By assumption, $m^*$ is an extreme point of $\mathcal{M}(G)$. Hence $m^* = m_k^*$ for all $k$ such that $\alpha_k > 0$. Therefore, all the $m_k$ lie in the preimage $(s^*)^{-1}(m^*) = \mathcal{M}(m^*)$. The point $m$ is an extreme point of $\mathcal{M}(m^*)$ and $m_k \in \mathcal{M}(m^*)$, and hence the equality $m = \sum \alpha_k m_k$ implies that there exists $k$ such that $m_k = m$ and $\alpha_k = 1$. This means that $m$ is an extreme point of the polyhedron $\mathcal{M}(Q)$.

This proves Proposition 6.

**Corollary 5.** *If* $m^* \in s^*(\mathcal{M}(Q))$ *is a simple flow, then* $m^*$ *is the image of some simple flow* $m$,

$$m^* = s^*m.$$

Consider successive subdivisions $C_1, C_2, C_3, \ldots$ such that the maximum diameter $d_1, d_2, d_3, \ldots$ of the cells in these divisions tends to zero. Here each covering $C_k$ consists of polyhedral cells intersecting in boundary discs. Let $C$ be one of the above coverings and let $C^* = \{M^*(i)\}$ be the division obtained from $C$ by attributing each boundary disc of neighbouring cells to one of the cells. The corresponding mapping $h(x) = \{i\colon x \in M^*(i)\}$ is uniquely defined. We assume without loss of generality that $C_1^*, C_2^*, C_3^*, \ldots$ are successive subdivisions of the original division $C_1^*$. Each invariant measure $\mu$ generates a flow on the symbolic image defined by

$$h^*\mu = m = \{m_{ij} = \mu(M^*(i) \cap f^{-1}(M^*(j)))\}.$$

So the flow $m^k = h_k^*\mu$ is defined on each symbolic image $G_k$. Moreover, the flows $m^k$ are consistent: $s^* m^{k+1} = m^k$.

Let $m_k$ be an arbitrary flow on $G_k$. We have the following alternatives: $m_k \notin s^*\mathcal{M}(G_{k+1})$ or $m_k \in s^*\mathcal{M}(G_{k+1})$. In the first case we say that the flow $m_k$ has no extension over $G_{k+1}$; in the second case, the flow $m_k$ is said to have an extension over $G_{k+1}$. According to Corollary 5, if a simple flow $m_k$ has an extension, then there exists simple flow $m_{k+1}$ such that

$$s^* m_{k+1} = m_k.$$

In other words, if a simple flow $m_k$ has an extension, then this extension can be chosen to be simple.

Suppose that on each symbolic image $G_k$ a flow $m^k$ is chosen and that these flows are consistent; that is,
$$s^* m^{k+1} = m^k.$$
For each $k$, consider the measure $\mu_k$ on $M$ defined by
$$\mu_k(A) = \sum_i m_i^k \frac{v(A \cap M(i))}{v(M(i))},$$
where $v$ is the Lebesgue measure.

**Theorem 9.** *Let $C_k$ be successive subdivisions such that the maximum diameter $d_k$ of the cells tends to 0. If $m^k$ is a consistent sequence of simple flows on the symbolic images $G_k$, then there exists an invariant measure*
$$\mu = \lim_{k \to \infty} \mu_k$$
*which is an ergodic measure.*

First we prove the following result.

**Proposition 7.** *Let $\mu_1 \neq \mu_2$ be measures. Then there exists a number $t^*$ such that $m_1^t = h_t^* \mu_1 \neq m_2^t = h_t^* \mu_2$ for all $t > t^*$, where $h_t^* \colon \mathscr{M}(f) \to \mathscr{M}(G_t)$ is the mapping of invariant measures into the space of flows on $G_t$.*

*Proof.* We proceed by *reductio ad absurdum*. Let $m_1^t = h_t^* \mu_1$ and $m_2^t = h_t^* \mu_2$ and suppose that there exists a sequence $t \to \infty$ such that $m_1^t = m_2^t$. By construction, each of the sequences $\{m_1^t\}$ and $\{m_2^t\}$ is consistent. For each $t$, let us construct measures $\mu_1^t$ and $\mu_2^t$ on $M$ by setting
$$\mu_{1,2}^t(A) = \sum_i (m_{1,2}^t)_i \frac{v(A \cap M(i))}{v(M(i))},$$
where $A$ is a measurable set, the $M(i)$ are cells in the division $C_t$ and $v$ is the Lebesgue measure. As a result, we get a sequence of measures $\{\mu_1^t\}$ and $\{\mu_2^t\}$ on the manifold $M$. By Theorem 1 in [19], there exist invariant measures
$$\mu_1 = \lim_{t \to \infty} \mu_1^t, \qquad \mu_2 = \lim_{t \to \infty} \mu_2^t.$$
By construction, $\mu_1^t = \mu_2^t$ for each $t$, and therefore $\mu_1 = \mu_2$. This contradiction completes the proof.

*Proof of Theorem* 9. By Theorem 1 in [19] the weak limit $\mu = \lim_{k \to \infty} \mu_k$ exists and the measure $\mu$ is invariant. Arguing by contradiction we show that $\mu$ is an ergodic measure. Let $\mu$ be written as
$$\mu = \alpha \mu_1 + (1 - \alpha) \mu_2,$$
where $\mu_1 \neq \mu_2$ and $0 < \alpha < 1$. Consider the mapping $h^* \colon \mathscr{M}(f) \to \mathscr{M}(G_k)$. We have
$$h^* \mu = \alpha h^* \mu_1 + (1 - \alpha) h^* \mu_2 = \alpha m_1^k + (1 - \alpha) m_2^k,$$

where $m_1^k = h^*\mu_1$ and $m_2^k = h^*\mu_2$. Since the sequence $m^k$ is consistent, we have

$$h^*\mu = m^k.$$

By Proposition 7, $m_1^k \neq m_2^k$ for all sufficiently large $k$. Hence

$$m^k = \alpha m_1^k + (1-\alpha)m_2^k,$$

where $0 < \alpha < 1$ and $m_1^k \neq m_2^k$. Therefore, the flows $m^k$ are not simple for large $k$. This contradiction completes the proof of the theorem.

**The extremal property of ergodic measures.** The set $\mathcal{M}(\Omega)$ of invariant measures concentrated on a component $\Omega$ of a chain recurrent set is a convex weakly compact set. The ergodic measures form the set of extreme points $\text{ext}\,\mathcal{M}(\Omega)$ of $\mathcal{M}(\Omega)$. By Choquet's theorem, any point in a convex compact set is an integral sum of its extreme points. We fix a continuous function $\varphi$ and consider the functional

$$\Phi(\mu) = \int_M \varphi\,d\mu$$

on the space of measures $\Phi\colon \mathcal{M}(\Omega) \to \mathbb{R}$. This functional is linear in $\mu$. As $\Phi$ is continuous in the weak topology, it maps the compact set $\mathcal{M}(\Omega)$ to a convex compact set in $\mathbb{R}$, namely the closed interval $[a,b]$, where

$$a = \min_{\mathcal{M}(\Omega)} \Phi(\mu) \quad \text{and} \quad b = \max_{\mathcal{M}(\Omega)} \Phi(\mu).$$

According to Theorem 2 in [25], $[a,b]$ is the spectrum of the averaging of the function $\varphi$ over periodic pseudotrajectories of the component $\Omega$.

**Theorem 10.** *Let $[a,b]$ be the spectrum of the averaging of a function $\varphi$ over the component $\Omega$. Then:*
- *the set*

$$\mathcal{M}(b) = \left\{\mu \in \mathcal{M}(\Omega)\colon \int_M \varphi\,d\mu = b\right\}$$

*is a nonempty convex compact subset of $\mathcal{M}(\Omega)$;*
- *the extreme points of $\mathcal{M}(b)$ are the ergodic measures $\mu^*$ such that $\int_M \varphi\,d\mu^* = b$;*
- *for any measure $\mu \in \mathcal{M}(b)$ there exists a division of the manifold $M$ (modulo nullsets) into invariant measurable subsets $\Omega_\alpha$, $\alpha \in A$, where $A$ is an abstract Lebesgue space with measure $\beta$ and each $\Omega_\alpha$ is equipped with an ergodic invariant measure $\mu_\alpha \in \text{ext}(\mathcal{M}(b))$, such that*

$$\int_M \theta\,d\mu = \int_A \left(\int_M \theta\,d\mu_\alpha\right) d\beta \tag{3.3}$$

*for any integrable function $\theta$.*
  *The same results also hold for $a = \min_{\mathcal{M}(\Omega)} \Phi(\mu)$.*

*Proof.* 1. According to Theorem 2 in [25], there exists a measure $\mu^*$ such that $b = \int_M \varphi \, d\mu^*$. Therefore, the set $\mathscr{M}(b)$ is nonempty. The functional $\Phi(\mu)$ is linear and continuous, and hence $\mathscr{M}(b) = \Phi^{-1}(b)$ is a convex compact set in the weak topology.

2. By Theorem 7,
$$b = \int_M \varphi \, d\mu^* = \int_A \left( \int_M \varphi \, d\mu_\alpha \right) d\beta, \qquad (3.4)$$
where $A$ is an abstract Lebesgue space with measure $\beta$, and moreover, there exists a division of $M$ (modulo nullsets) into invariant subsets $\Omega_\alpha$, $\alpha \in A$, and each $\Omega_\alpha$ is equipped with an ergodic measure $\mu_\alpha$. We set
$$\psi(\alpha) = \int_M \varphi \, d\mu_\alpha.$$
Hence $\psi(\alpha)$ is a measurable function and $b = \int_A \psi(\alpha) \, d\beta$. We fix $\varepsilon > 0$ and construct the $\beta$-measurable sets $A_1 = \psi^{-1}([a, b-\varepsilon])$ and $A_2 = \psi^{-1}((b-\varepsilon, b])$. The set $A$ can be written as $A = A_1 + A_2$. Arguing by contradiction, we show that the measure $\beta(A_1)$ is zero. Assume on the contrary that the $\beta$-measure of $A_1$ is $\delta > 0$. Then $\beta(A_2) = 1 - \delta$. We have
$$b = \int_A \psi(\alpha) \, d\beta = \int_{A_1} \psi(\alpha) \, d\beta + \int_{A_2} \psi(\alpha) \, d\beta \leqslant (b-\varepsilon)\delta + b(1-\delta) = b - \varepsilon\delta.$$
If $\delta > 0$, then we obtain $b \leqslant b - \varepsilon\delta$, which is a contradiction.

Hence the $\beta$-measure of the set $A_2 = \psi^{-1}((b-\varepsilon, b])$ is 1 for any $\varepsilon > 0$. This gives us a monotone sequence of nested sets $\{\psi^{-1}((b-\varepsilon, b]), \varepsilon > 0\}$ such that
$$\bigcap_{\varepsilon > 0} \psi^{-1}((b-\varepsilon, b]) = \psi^{-1}(b).$$
By monotonicity,
$$\beta(\psi^{-1}(b)) = \beta\left(\bigcap_{\varepsilon > 0} \psi^{-1}((b-\varepsilon, b])\right) = \lim_{\varepsilon \to \infty} \beta(\psi^{-1}((b-\varepsilon, b])) = 1.$$
The Lebesgue space $A$ is defined up to nullsets with respect to the $\beta$-measure. Hence we can assume that $A$ coincides with
$$\psi^{-1}(b) = \left\{ \alpha \colon \mu_\alpha \in \operatorname{ext}(\mathscr{M}(\Omega)), \int_M \varphi \, d\mu_\alpha = b \right\}.$$
By the above, the family of ergodic measures
$$\left\{ \mu_\alpha \in \operatorname{ext}(\mathscr{M}(\Omega)) \colon \int_M \varphi \, d\mu_\alpha = b \right\}$$
forms the set of extreme points of the convex set $\mathscr{M}(b)$.

3. The third assertion of the theorem follows from assertion 2 of this theorem and Theorem 7.

Theorem 10 is proved.

The numbers $a$ and $b$ and the measures $\mu_a$, $\mu_b$ can be evaluated to any accuracy using the symbolic image as follow. Let $G$ be a symbolic image for the covering $C = \{M(i)\}$. The function $\varphi$, which generates the framing of the symbolic image, associates the number $b(i) = \varphi(x_i)$, $x_i \in M(i)$, with each vertex $i$. This limit set of the averagings over all periodic paths of the component $H(\Omega)$ is the spectrum of the framing. Theorem 1 in [25] asserts that the spectrum of the framing $[\alpha, \beta]$ converges to $[a, b]$ as the diameter of the covering $C$ goes to zero. Moreover, the numbers $\alpha$ and $\beta$ are the averagings over the simple cycles $\omega_\alpha$ and $\omega_\beta$, respectively. The cycles $\omega_\alpha$ and $\omega_\beta$ define simple flows on $G$, which approximate the ergodic measures $\mu_a$ and $\mu_b$. The cycles $\omega_\alpha$ and $\omega_\beta$ can be found numerically using the approach in [26]–[28]. This means that the ergodic measures $\mu_a$ and $\mu_b$ can be approximated numerically.

**Georgy S. Osipenko**
Sevastopol Branch
of Lomonosov Moscow State University
*E-mail*: george.osipenko@mail.ru